\documentclass[11pt]{amsart}
\usepackage{amssymb}
\usepackage{amsmath,latexsym}

\newenvironment{prf1}{{\it Proof of Theorem  \ref{mainthm}.}}
{\hfill $\Box$ \medskip}
\newenvironment{prf2}{{\it Proof of Corollary  \ref{CPWcrl}.}}
{\hfill $\Box$ \medskip}
\newenvironment{prf3}{{\it Proof of Corollary \ref{maincor}.}}
{\hfill $\Box$ \medskip}

\newtheorem{thm}{Theorem}[section]

\newtheorem{dfn}[thm]{Definition}
\newtheorem{lmm}[thm]{Lemma}
\newtheorem{crl}[thm]{Corollary}

\newtheorem{rms}[thm]{Remarks}

\def\Irr{\mathrm {Irr}}
\def\IBr{\mathrm {IBr}}
\def\CF{\mathrm {CF}}
\def\Res{\mathrm {Res}}

\def\Br{{\rm Br}}

\def\Oo{{\mathcal O}}

\begin{document}

\title{On isotypies  between  Galois conjugate blocks}
\author{Radha Kessar}

\address {Institute of Mathematics, University of Aberdeen, 
Fraser Noble Building, King's College, Aberdeen AB24 3UE, U.K. }

\begin{abstract} We show that between any pair of Galois conjugate  blocks of 
a  finite group, there is an isotypy with all signs positive.
\end{abstract}
%\e-mail{kessar@maths.abdn.ac.uk}
\maketitle

\section{Introduction}

Let $p$ be a prime number, let $k$ be an  algebraic closure
of the field of $p$ elements  and  let $G$ be a finite group. 
Let $$\sigma:  k \rightarrow k, \ \  (\lambda \to  \lambda^ p,  
\ \ \lambda \in k )$$ be the  
Frobenius homomorphism of  $k$ and
we write also
$\sigma: kG \rightarrow kG$ for the   ring automorphism  
induced by  $\sigma $. This is defined by 
\[ \sigma\Bigl(\sum_{g\in G}\alpha_g g\Bigr)= \sum _{g\in 
G}\alpha_g^p g. \] 
The map $\sigma$ is a  ring automorphism but
not a $k$-algebra  automorphism. Since   $\sigma $ is a ring automorphism, 
$\sigma $  induces a permutation of  the  blocks  of $kG$.  Here 
by a block of a ring $A$, we mean a primitive idempotent of the center of   
$A$.  Blocks $b $ and $ c$ of  $kG$ are said to be  {\em Galois 
conjugate}  if $ c= \sigma^n(b)$ for some natural number $n$.

Let $(K, \Oo , k)$  be a $p$-modular system  and let 
$\,^- : \Oo \to   k $  be the canonical quotient mapping. Let   
$\,^- : \Oo G \to   kG $   be the induced $\Oo$-algebra homomorphism of the 
group algebras.   In particular, $ b \to \bar b $  induces a bijection 
between the  set of  blocks  of $\Oo G$ and the  set of  blocks  of $kG$.
Blocks $b $ and $ c$ of  $\Oo G $ are said to be  {\em Galois 
conjugate}  if $\bar b $ and $\bar c$ are Galois conjugate   blocks of $kG$.

Our main result is the following.

\begin{thm}   \label{mainthm}     With the notation above, suppose that 
$K$ contains a primitive  $|G|$-th  root of unity. 
Then any two  Galois conjugate blocks of $\Oo G$ are isotypic.
\end{thm}

\iffalse
For any  block $B$ of $kG$, $\sigma $ defines an isomorphism of $k$-algebras 
between $B^{(p)}$ and $\sigma(B)$
as $k$-algebras.   
Here $B^{(p)}$  is the Frobenius twist of $B$ defined  as 
follows:
The underlying ring of $B^{(p)}$ is the same as that of   $B$, but
we endow it with a new action of the scalars in $k$ via the Frobenius
map on $k$: for $\lambda\in k$ and $x\in A$, the new action is given by
$\lambda\cdot x=\lambda^{\frac{1}{p}} x$. 
\fi

As a corollary, we obtain a  new proof  of   the  following result
of Cliff, Plesken and Weiss \cite{CliffPleskenWeiss}.  

\begin{crl}\label{CPWcrl}  
Let $G$ be a finite group and $b$ be a block of $kG$. Then the center
$Z(kGb)$   of $kGb$ has an ${\mathbb F}_p$-form. 
\end{crl}

We recall that a finite dimensional  $k$-algebra $A$ is said to have an  
${\mathbb F}_p$-form if   there exists an 
$ {\mathbb F}_p$-algebra  $A_0$ such that   
$ A \cong k \otimes_{{\mathbb F}_p} A_0 $ as  $k$-algebra.

The notion of isotypic blocks  (see   Definition \ref{defnisot})  is due to 
Michel Brou\'e \cite[Definition 4.3]{broue-type}.  
Isotypies   between blocks   are  interesting   as they are  often 
the character theoretic shadow of  a    derived  or Morita equivalence 
between the   $k$-linear module categories  of the  corresponding 
block algebras.  In \cite{Benson/Kessar:2007}, it was shown  that 
there exist pairs of  Galois conjugate blocks 
which are not  derived equivalent as $k$-algebras.
Thus the blocks studied in \cite{Benson/Kessar:2007} provide examples
of pairs of blocks  which are isotypic  and isomorphic as rings  
but not derived equivalent as $k$-algebras. 
However, it is conjectured that  the   number of Morita equivalence classes
of algebras in any set of Galois conjugate blocks is bounded by  a number
which   depends only on  the  defect of the  blocks and  
which is independent of $G$. This conjecture is   related to  the  
Donovan finiteness conjectures in block theory \cite{Ke1}.

For a   non-negative   integer  $d$ and  a finite dimensional commutative 
$k$-algebra $A$, 
we will say that  $A$  occurs as the center of a $d$-block if there exists a 
finite group $H$ and  a block $c $ of $kH$ with defect $d$ such that    
$A \cong Z(kHc)$ as  $k$-algebras.
Combining  Corollary \ref{CPWcrl} with a   theorem  of 
Brauer  and Feit \cite{brauerfeit} gives  a proof  of 
the following  finiteness  result. The result itself is known to  
experts, but as far as I am aware has not appeared  before in the 
literature.

\begin{crl}\label{maincor}  Let  $d$ be a  non-negative integer and set  
$m:= \frac{1}{4} p^{2d}+ 1 $. Up to isomorphism there are  at most 
$p^{m^3}$ $k$-algebras  that occur  as centers  of  $d$-blocks.
\end{crl}

The paper is divided into three sections.  In Section 2, 
we   set up notation and recall  the definitions of perfect isometries 
and isotypies. Section 3 contains the proofs of Theorem \ref{mainthm},
Corollary \ref{CPWcrl}  and 
Corollary \ref{maincor}.

\section{Notation and definitions} 
Throughout, $G$ is a finite group,  $k$ is algebraically closed of  
characteristic $p$  and 
$(K, \Oo , k)$  is a $p$-modular system such that   $K$ contains a primitive 
$|G|$-th root of unity.  For an element $ g $ of  $G$, we  will denote 
by $g_p $ the 
$p$-part of 
$g$ and by $g_{p'}$-the $p'$-part of $G$.  Denote by $G_p $ the 
set of elements of $G$ of  $p$-power order and by $G_{p'}$ 
the set of elements of $ G$ of $p'$-order.  For a natural number $n$, let 
$n_p$ denote the $p$-part  $n$ and $n_{p'}$ denote the  $p'$-part of $n$. 

\subsection {Generalized decomposition maps} 

Let $b$ be a  central idempotent  of ${\mathcal O}G$. 
Denote by  $\CF(G, K)$  the $K$-space of $K$-valued  class  functions  on 
$G$.  Identifying   $K$-valued  class functions on $G$ 
with their canonical $K$-linear extensions to $K$-valued  functions on $KG$,
denote     by  $\CF(G, b, K)$   
the $K$-subspace of class functions $\phi $  in $\CF(G, K)$ such that 
$\phi(gb) = \phi(g)$ for all $g \in G$.  Denote by  
$\CF_{p'}(G, K)$  (respectively  $\CF_{p'}(G, b, K)$)
the $K$-subspaces of $\CF(G, K)$
(respectively $\CF(G, b, K))$ of class functions  
which vanish on  the $p$-singular 
classes of $G$;  
$\Irr (G, K)  $  (respectively     
$\Irr(G, b, K)$)   the subset of   $\CF(G, K)$ (respectively  $\CF(G, b,K)$)   
of irreducible $K$-valued characters of $G$;   by  $\IBr(G, K)$  
(respectively $\IBr(G,b, K)$) the set of irreducible Brauer   characters of 
$G$ (respectively 
irreducible  Brauer characters of $G$  in $b$).  For $\phi \in \IBr(G, K)$, let 
$\phi_0 \in \CF_{p'}(G, K)$  be defined by  $\phi_0(y) =\phi(y) $  if   
$y \in G_{p'}$.   Denote by 
$\IBr(G, K)_0 $ the set  of class functions $\phi_0$, for 
$\phi \in \IBr(G, K)$ and by 
$\IBr(G,b, K)_0 $ the  subset  of $\IBr(G, K)_0 $ 
consisting of those $\phi_0$ for which   $\phi \in \IBr(G,b, K)$

\iffalse
and  let $ \alpha \in \CF(G, K)$. Denote  by 
$\alpha(x-) $ the element of $\CF_{p'}(C_G(x), K) $ defined by  
$\alpha(x-)(y) =  \alpha(xy)$  for  $y \in C_G(x)_{p'}$.
\fi

Let $x$  be a  $p$-element of $G$.  
The  generalized decomposition map
$$d_G^x: \CF(G, K)  \to \CF_{p'}(C_G(x), K) $$ is   defined by 
$d_G^x(\alpha)(y) = \alpha(xy) $ for $\alpha \in   \CF(G, K)   $, $ y \in    
 C_G(x)_{p'}$.  If   $\chi \in   \Irr(G, K)$,  $x  \in G_p$, then  
$$d_G^x(\chi)       = \sum_{\phi \in \IBr(C_G(x), K)} \delta_{(G,\chi)} 
^{(x,\phi)}\phi_0 $$ 
for uniquely determined elements $\delta_{(G,\chi)} ^{(x,\phi)}  $ of $K$; 
$\delta_{(G,\chi)} ^{(x,\phi)}$ is called 
the generalized decomposition number  associated to  the triple 
$(\chi, \phi, x ) $.

\iffalse
If $e$ is a central idempotent of ${\mathcal O}C_G(x)$, and 
$ \alpha \in \CF(G, K)$, denote by $\alpha(xe-) $ the element of 
$\CF_{p'}(C_G(x), K) $ defined by  
$\alpha(xe-)(y) =  \alpha(xye)$  for  $y \in C_G(x)_{p'}$.
\fi

Let $e$ be  a central idempotent of ${\mathcal O}C_G(x)$.  The map
$$d_G^{(x,e)}: \CF(G, K)  \to \CF_{p'}(C_G(x),e, K) $$  is defined by
$d_G^{(x,e)}(\alpha)(y) = \alpha(xey) $ for   
$\alpha \in \CF(G,K)$, $y\in    C_G(x)_{p'}$.

\subsection{Local structure of blocks}  Let $R$ denote one of the 
rings $\Oo $ or $k$. For  $Q$  a $p$-subgroup of $G$, 
 let $ (RG)^Q $  denote    the    $R $-subalgebra of 
$RG $ consisting  of  the elements of  $ RG$ 
which are fixed by the conjugation action of $Q$.
The Brauer homomorphism  $\Br_Q: (RG)^Q  \to kC_G(Q)$  is the map 
defined by 
$$\Br_Q (\sum_{g\in g}\alpha_g g ) = \sum_{g\in kC_G(Q) }\bar{\alpha}_gg. $$
Here if $R=k$, then $\bar {\alpha}_g $ is to be interpreted as $\alpha_g $.
If $b$  is a block of ${\mathcal \Oo}G$,  then a $\bar b$-Brauer pair  
is a pair 
$(Q, \bar e)$, where  $Q$ is a $p$-subgroup of $G$  and $e $ is a block of 
${\mathcal O}C_G(Q)$ such that $\Br_Q(b) \bar e = \bar e $.  
Let  $(P, \bar e_P)$  be  a maximal  $b$-Brauer  pair under the 
Alperin-Brou\'e inclusion of Brauer pairs \cite[Definition 1.1]{alperinbroue}, 
and for each  subgroup $Q$ of $P$, let $e_Q$ be the unique block of 
${\mathcal O}C_G(Q)$   such that $(Q,\bar e_Q) \leq (P, \bar e_P)$.

Let  ${\mathcal F}_{(P, \bar e_P)}(G,\bar b)$ denote the  category  
whose objects are the  subgroups of $P$ and whose morphisms are defined as 
follows: For $Q, R$ subgroups of $P$,   the set  of 
${\mathcal F}_{ (P, \bar e_P)}(G,\bar b)$  morphisms  from $Q$ to $R$  
is the set 
of those    group homomorphisms $ \varphi: Q \to R$ for which there exists 
an element 
$g$ of $G$ such that  $\varphi(x) =gxg^{-1}$ for $x \in Q$ and   such that  
$ (\,^gP, \,^g\bar e_Q )\leq (R, \bar e_R)$;  composition of morphisms 
 is  the  usual  composition of  functions.

\subsection{Perfect isometries and isotypies} Let    
$H$ be a finite group such that $K$ contains a primitive $|H|$-th root of 
unity. Let $b$ be a  central idempotent of ${\mathcal O}G$ and 
$c$ a   central idempotent  of ${\mathcal O}H$.

\begin{dfn} A perfect isometry between  $b$ and $c$ is a $K$-linear map
$$I : \CF(G,b, K) \to \CF(H,c, K)   $$
such that the following holds.
For each $\chi  \in \Irr(G,b,K)$, there exists an 
$\epsilon_{\chi} \in \{ \pm 1\} $ such that   
 the map 
$ \chi \to \epsilon_{\chi}I(\chi) $ is a bijection between   $\Irr(G,b,K)$ and 
 $\Irr(H,c,K)$ and such  that setting 
$$\mu := \sum_{ \chi \in \Irr(G,b,K)}  \chi \times I(\chi), $$ 
the class function on $G \times H$ which sends an element $(x,y)$ of  
$ G\times H$ 
to the element   $\sum_{ \chi \in \Irr(G,b,K)}  \chi(x) I(\chi)(y)$   of 
${\mathcal O}$   the following holds:

(a) For each $x\in G$, $y\in H$,  
$\frac{\mu(x,y)} {|C_G(x)|} \in  {\mathcal O}$.

(b)   If $x \in G$, $y\in H$ are such that exactly one of $x $ and $y$ is  
$p$-singular, then $ \mu(x,y)=0 $.
\end{dfn}

If  $I$  as  in the above  definition  is a perfect isometry between 
$b$ and $c$, then  $I$ induces by restriction
a map
$$I_{p'} : \CF_{p'}(G,b, K) \to \CF_{p'}(H,c, K).   $$

\begin{dfn}\label{defnisot}   Let  $b$  be a block  of ${\mathcal O}G $ and 
 $c$  a block  of ${\mathcal O}H $. Then  $b$ and $c$  are isotypic   
if the following conditions hold:

(a) There exists a  $p$-group $P$ and inclusions 
$ P \hookrightarrow G$,  $P\hookrightarrow H$ such that identifying  $P$ 
with  its  image in $G$ and in  $H$,  there exists a block $e_P$ of 
${\mathcal O}C_G(P)$ such that 
$(P, \bar e_P) $ is a maximal $\bar b$-Brauer pair and  
a block 
$f_P$ of ${\mathcal O}C_H(P)$ such that $(P, \bar f_P)$ is a maximal   
$\bar c$-Brauer pair and   such that 
$${\mathcal F}_{(P,\bar e_P)}(G,\bar b)  = 
{\mathcal F}_{(P, \bar f_P)}(H,\bar c). $$

(b) For each cyclic subgroup $Q$ of $P$, there exists a perfect isometry
$$I^Q: \CF(C_G(Q), e_Q,  K) \to \CF(C_H(Q), f_Q,  K) $$ where $e_Q$ 
(respectively $f_Q$) is the unique block of ${\mathcal O}C_G(Q)$ 
(respectively ${\mathcal O}C_H(Q)$)  with 
$ (Q, \bar e_Q) \leq (P,\bar e_P) $ (respectively 
$ (Q, \bar f_Q) \leq (P,\bar f_P)$) such that 
for every generator $ x$ of $Q$, we have
$$ I^Q_{p'} \circ d_G^{(x, e_Q)}  =  d_H^{(x, f_Q)}\circ  I^{\{1\}} . $$  
\end{dfn}

\section{Proofs.}  
Keep the notation  and hypothesis of Theorem \ref{mainthm}.    
 In addition, let  $W(k)$ be the unique 
absolutely non-ramified complete   discrete  valuation ring   having $k$ as 
residue field, and  identify $W(k)$ with its image   under the canonical  
injective homomorphism   $W(k) \rightarrow {\mathcal O}$ 
(see \cite [Chapter 2 \S 5, Theorem 3 and Theorem 4]{Se}).     
There is a unique ring  automorphism  
$\sigma_{W(k)} :  W(k) \rightarrow W(k)$  such that    
$ \overline{\sigma_{W(k)}(\eta)}= 
\sigma(\bar \eta)$ for all $\eta \in W(k)$. Further,  note that for any 
$p'$-root of unity  $\eta$ in $K$, $\eta \in W(k)$ and 
$\sigma_{W(k)}(\eta)= \eta^p$.

Let $K_0 $  be the  algebraic closure of  ${\mathbb Q}$ in  $K$.   
Choose a field automorphism  $\sigma_{K_0}   : K_0 \to K_0 $  such that    
if $ \eta \in K$ is any $|G|$-th  root of unity, then 
$\sigma_{K_0} (\eta) = \eta^p$   if  the order of $\eta $ is relatively prime  
to $p$  and $\sigma_{K_0} (\eta) =\eta $  if $\eta $ is  a power 
of $p$.
Then, $ K_0 \cap W(k)$ contains   a  primitive  $|G|_{p'}$-root of unity   and 
$\sigma_{K_0}$ and $ \sigma_{W(k)}$   coincide on    
${\mathbb Q}[\eta ]\cap  W(k)$ for any $|G|_{p'}$-root of unity $\eta$.  
Note that we are not claiming that  $\sigma_{0}$  extends to 
an automorphism of $\Oo $  which agrees with $\sigma_{K_0}$ on restriction to 
${\mathbb Q}[\eta ]\cap  \Oo$ for any $|G|$-th root of unity $\eta$.

Denote by  $\sigma_{W(k)}$ (respectively $\sigma_{K_0}$)
the natural  extension  of  $\sigma_{W(k)}$ (respectively $\sigma_{K_0}$) to  
$W(k)G$ (respectively $K_0G $).

Recall from \cite[Chapter 3, Theorem 6.22 (ii)]{Nagao-Tsushima} that    
any block of   $\Oo H$, for $H$ a finite group  is an 
$\Oo $-linear combination of $p'$-elements of $H$. 

\begin{lmm} \label{blockmoves}  Let   $\eta $ be a  primitive $|G|_{p'} $  
root of unity   
in $K$ and let  $H \leq G$. Let
$$ c = \sum_{g \in H_{p'}} \alpha_g g,   \alpha_g \in {\mathcal O}$$   
 be   a   block  of ${\mathcal O}H$. Then, 

(i) $\alpha_g   \in   {\mathbb  Q}[\eta]  \cap W(k)$ for all $g\in H_{p'}$.

(ii)  $\sigma_{W(k)}(c)$ is a  block  of $\Oo H$,  $\sigma_{W(k)}(c)=
\sigma_{K_0}(c)$  and 
$$\overline {\sigma_{W(k)}(c)} = \sigma(\bar c). $$
\end{lmm}

\begin{demo}(i) By idempotent lifting, the   canonical quotient map  
$\,^-: \Oo \to k $ induces  a bijection between the 
set of central idempotents of $\Oo H$ and $kH$. Similarly, the  
restriction of $\,^- $   to $W(k)$ induces a bijection between the 
set of central idempotents of $W(k) H$ and $kH$. 
Since a central idempotent of $W(k)H$ is  a central idempotent of $\Oo H$, 
it follows that  $\Oo H$ and $W(k)H$ have the same central idempotents. 
In particular, $c  \in W(k) H$. So,   $\alpha_g  \in W(k)$ for all  
$ g\in H_{p'}$.  On the other hand, by \cite[Chapter 3, Theorem  6.22]
{Nagao-Tsushima},  
for  $g\in H_{p'}$,  
$\alpha_g $ is a ${\mathbb Q}$-linear combination of $|g|$-th roots of unity  
whence 
$\alpha_g \in {\mathbb Q}[\eta] $.  This proves (i).
 
(ii) As shown above, the set of blocks  of $\Oo H$  is the same as the set of 
blocks of 
$W(k)H$  and   $\sigma_{W(k)} $ is an automorphism of 
$W(k)H$.   This proves the first assertion. The others are immediate from (i). 
\end{demo}

\

\begin{dfn}  \label{mapofiso} For $H$   a subgroup  $G$, let       
$${I}^H : \CF(H, K) \to  \CF(H,K)$$ 
denote  the $K$-linear map  defined by
$${I}^H(\phi)(x) = \phi(x_p x_{p'}^p), 
\ \ \text{ for } \phi \in \CF(H,K), \ x  \in G  . $$
\end{dfn}

\

If an element  $\phi  \in \CF(H, K)$
takes values in $K_0 $,   denote by   
$\sigma(\phi)$  the  element of $\CF(H, K)$   which  sends  $g \in H $  to 
$\sigma_{K_0} (\phi(g)) $.  Similarly, if $ \phi $ takes values in $W(k)$, 
then   $\sigma_{W(k)}(\phi) $ will denote   the     element of $\CF(H, K)$  
which  sends  $g \in H $  to $\sigma_{W(k)} (\phi(g)) $.   
We use the same conventions  on  $K$-valued class  functions defined  on 
$H_{p'}$.

\begin{lmm} \label{isomoves}  Let $H$ be a subgroup of $G$ and let 
$c$ be a block of $H$.

(i)   For any  $\chi \in \Irr (H, K)$,  
$${I}^H(\chi) = \sigma_{K_0}(\chi ) \in   \Irr(H,K). $$
 
(ii)  The map  $$\chi   \to  \sigma_{K_0}(\chi)    $$  
is a bijection  from  
$ \Irr (H,c, K) $    to  $ \Irr (H, \sigma_{K_0}(c), K)  $.

(iii) The  restriction of ${I}^H $ to $\CF(H,c,K)$ 
induces a perfect isometry between $c$ and $\sigma_{K_0}(c)$.

(iv)  For any  $\phi \in \IBr(H,K)$ and any $\chi \in \Irr(G,K)$, 
$${I}^H(\phi_0)  =  \sigma_{K_0}(\phi)_0  \in  \IBr (H, K)_0 $$  
and $$\delta_{ (H,  \sigma_{K_0}(\chi))}^{ (\{1\},  \sigma_{K_0}(\phi))} = 
\delta_{(H,  \chi)}^{ (\{1\},  \phi)} . $$
\end{lmm}

\begin{demo}  (i)  
Let $\chi \in \Irr(H,K)$, and let $\rho: H \to GL_n(K_0)$ be a 
representation 
affording $\chi $.    Such a $\rho $ exists by   Brauer's  splitting 
field theorem 
(see for example \cite[Chapter 3, Theorem 4.11]{Nagao-Tsushima}).  
Denoting  also by 
$\sigma_{K_0}$ the automorphism (as abstract group) of $GL_n(K_0)$ 
induced by $\sigma_{K_0}$, 
it follows that  $\sigma_{K_0} \circ \rho $ is an irreducible representation of 
$H$  with character $\sigma_{K_0}(\chi )$.  
We  show that ${I}^H (\chi) = \sigma(\chi)$.
Let $ h =h_ph_{p'}  \in H$.   Since $K_0$ contains all   the eigen values of 
$\rho (h)$,  by   replacing $\rho $
with  an equivalent representation if necessary, we may assume that 
$\rho(h) $  is a diagonal matrix.   Since  $h_p$ and $h_{p'}$ are powers of 
$h$, 
it follows that $\rho(h_p) $   is a diagonal matrix 
$ diag(\zeta_1, \cdots, \zeta_n) $, 
$\rho(h_{p'}) $   is a diagonal matrix $diag(\eta_1, \cdots, \eta_n) $, 
where  each $\zeta_i $ is a  $|H|_p$ root of unity and each $\eta_j $ is an 
$|H|_{p'}$ root of unity,
$\rho (h)$is the diagonal matrix
$ diag(\zeta_1\eta_1, \cdots, \zeta_n\eta_n  ) $ and 
$\rho (h_ph_{p'}^p)$ is the diagonal matrix
$ diag(\zeta_1\eta_1^p, \cdots, \zeta_n\eta_n^p  ) $. 
Thus, 
$${I}^H(\chi)(h)=  \sum_{1\leq i \leq n } \zeta_i\eta_i^p = 
\sigma_{K_0} (\sum_{1\leq i \leq n } \zeta_i\eta_i) = \sigma(\chi) (h). $$

(ii)  Let   $ \chi \in \Irr (H, K)$  and let 
$$e_{\chi} = \frac{\chi(1)}{|H|}\sum_{h \in H } \chi (h)  h^{-1} $$
be the primitive central idempotent  of $KH$ corresponding to $\chi $. 
Then, 
$$e_{\sigma_{K_0} (\chi)}= \frac{\chi(1)}{|H|}\sum_{h \in H } 
\sigma_{K_0} (\chi (h)) h^{-1}  
=  \sigma_{K_0} (e_{\chi}).$$
From this it is immediate that if $ \chi \in \Irr(H,c, K)$
then $\sigma_{K_0}(\chi) \in \Irr(H,\sigma_{K_0}(c), K)$.

(iii)  By  (i)  and (ii) it suffices to prove that the   function
$$\mu:= \sum_{\chi \in \Irr(H,c,K)} \chi \times {I}^H (\chi), $$
satisfies conditions  (a) and (b) of Definition 1.1.

Set $$\iota := \sum_{\chi \in \Irr(H, K, c)}\chi \times  \chi  $$
and let $x, y \in H$. 

Then   $$\mu(x, y) = \iota (x,  y_p y_{p'}^p). $$ Since the 
identity map on $\CF(H, c , K)$ is  a perfect isometry with  $\iota $ 
the corresponding (virtual) character  of $H \times H$,  we see that
$\mu(x, y)$ is divisible in ${\mathcal O}$ by both $|C_G(x)|$ and by 
$|C_G( y_p y_{p'}^p  )| $ and that 
$\mu(x, y)$ is  $0$ if exactly one of $x$ and 
$y_p y_{p'}^p $ is $p$-singular. But, clearly  $C_G( y_p y_{p'}^p)= C_G(y)$ 
and 
$y_p y_{p'}^p $ is $p$-singular if and only if $y$ is $p$-singular. 
This proves (iii).

(iv) Let  $\phi \in \IBr(H, K) $ and let  $\tau: H \to GL_n(k)  $  
be an irreducible representation of  
$H$   affording the  Brauer character $\phi $.
Then $\sigma \circ  \tau $ is an irreducible representation of $H$, 
where again we denote by $\sigma $ the induced automorphism  of 
$GL_n(k)$.   Let $\phi'$ be the Brauer character associated  to 
$\sigma \circ \tau $ and let $ h \in   H_{p'}$.  Let  $\Lambda $ be the 
multiset  of   
eigen values of  $\tau(h)$.
The multiset of eigen values  of $\sigma \circ \tau (h) $ is 
$\{ \lambda^p \, :  \lambda \in \Lambda \}$, hence $\phi'(h) =  
\sigma_{W(k)}(h)$.  
Further, since  $\{ \lambda^p \, :  \lambda \in \Lambda \}$
is also the multiset of eigen values of $\tau(h^p)$,     
$ {I}^H(\phi_0) (h) = \phi(h^p) =\sigma_{W(k)}(h)$  for all 
$ h \in H_{p'}$.    Thus, 
$${I}^H(\phi_0) =  \sigma_{W(k)}( \phi)_0 =\phi'  \in  \IBr(H,K).$$
Since $\phi $ takes values in  ${\mathbb Z}[\eta ]$ , with  $\eta $  
a primitive 
$|H|_{p'}$-root of unity,   $\sigma_{W(k)} (\phi )_0=\sigma_{K_0}(\phi)_0 $  
whence ${I}^H(\phi_0) =  \sigma_{K_0}(\phi)_0   \in  \IBr(H,K)$.

The compatibility of decomposition numbers is clear from the  
fact that for any $\chi   \in \Irr (H,K)$,
$$ \Res|_{H_{p'}} \chi = \sum_{\phi \in \IBr(H,K)}
\delta_{H,  \chi}^{{1},  \phi}  \phi $$
and that  $\delta_{H,  \chi}^{{1},  \phi} \in {\mathbb Z}$ for all 
$\phi \in \IBr(H,K)$.

\end{demo}

\

\begin{prf1} Let $b$ be a block of ${\mathcal O}G $.   
By Lemma \ref{blockmoves}  it  suffices  to show that   
$b$ and $\sigma_{K_0}(b)$  are isotypic.     Let $P$ a $p$-subgroup of 
$G$  and   $e_P$ be a block of ${\mathcal O}C_G(P)$ 
such that $(P, \bar e_P)$ is a maximal $\bar b$-Brauer pair.   
For each $ Q\leq P$, let $e_Q $ be the unique block of 
${\mathcal O}C_G(Q)$ such that  $(Q, \bar e_Q) \leq  (P, \bar e_P)$.       For any   $p$-subgroup $Q$ of $G$  and any $a \in  (kG)^Q $, $\sigma(\Br_P(a)) =  \Br_P(\sigma (a) )$.  So, the map  $( Q, f) \to  (Q, \sigma(f)) $   is an isomorphism   from  the $G$-poset of 
$\bar b $-Brauer pairs   to the $G$-poset of $\sigma (\bar b)$-Brauer pairs.    By Lemma  \ref{blockmoves},  $\sigma (\bar b)=  \overline{\sigma_{K_0}(b)}$  and $\sigma( \bar f) =  \overline{\sigma_{K_0}(f)} $  for any block $f$ of $kC_G(Q)$, $Q$ a $p$-subgroup of   $G$.
Thus, $(P,\overline{ \sigma_{K_0}(e_P)}) $ is a maximal $\overline{\sigma_{K_0}( b)}$-Brauer pair;  for every  subgroup $Q $ of $P$, $\sigma_{K_0}( e_Q) $ 
is the unique block
of $\Oo C_G(Q)$ such that  $(Q, \overline{\sigma_{K_0}( e_Q)} ) \leq (P,  \overline{ \sigma_{K_0}(e_P)}) $, and 
${\mathcal F}_{(P,       \overline{\sigma_{K_0}(e_P) }   ) }(G,  \overline{\sigma_{K_0}(b)})  = {\mathcal F}_{(P, \bar e_P)}(G,\bar b)$.

We will  use the  maps ${I }^{C_G(Q)}$  of Definition \ref{mapofiso} 
to produce  an isotypy between  $b$ and  $\sigma_{K_0}(b)$.   For  
$Q \leq P$, let  $I^Q$ be the restriction of ${I}^{C_G(Q)}$  to 
$\CF(C_G(Q), e_Q, K)$.  We will show that   the   maps $I^Q$, 
as $Q$ runs over the cyclic subgroups of $P$ defines  an isotypy between  
$b$ and $c$.   So, let $Q \leq P$ be a  cyclic group. 
By Lemma \ref{isomoves}, 
$$I^Q  : \CF(C_G(Q), e_Q, K)  \to  \CF(C_G(Q), \sigma_{K_0} (e_Q), K) $$
is  a perfect  isometry.  It remains only   to  check the compatibility 
condition (b) of Definition \ref{defnisot}.  

Set $H=C_G(Q)$ and let   $Q= \langle x \rangle $.
We claim that 
$$\delta_{(G,  \chi)}^{ (x,  \phi)} 
=  \delta_{ (G,  \sigma_{K_0}(\chi))}^{ (x,  \sigma_{K_0}(\phi))} $$ 
for all 
$\chi \in \Irr(G,  K)$ and  all $\phi \in \IBr(H, K)$.  

Indeed, for   
$\tau  \in \Irr(H, K)$, let  $\zeta_{(\tau, x)} \in K $  
be such that  $\rho(x) $ is the scalar matrix  
$(\zeta_{(\tau, x)}, \cdots, \zeta_{(\tau, x)})$ in any representation of 
$H$ affording $\tau $.     Since the order of $\zeta_{(\tau, x)} $ 
is a divisor of $|G|_p $, 
$$\zeta_{(\sigma_{K_0}(\tau), x)}=  \sigma_{K_0}(\zeta_{(\tau, x)})
= \zeta_{(\tau, x)}. $$ 
Also, for $\chi \in   \Irr(G,K)$,  $\phi \in  \IBr(H, K)$, 
$$\delta_{(G,  \chi)}^{ (x,  \phi)}=  \sum_{ \tau \in \Irr(H, K)} \langle 
\Res|_H\chi, \tau \rangle \zeta_{(\tau, x)}\delta_{(H,  \tau)}^{ (\{1\},  \phi)}.$$

The claim follows  since for  all 
$\tau \in \Irr(H,K)$ and all $\phi \in  \IBr(H,K)$, 
$\zeta_{(\sigma_{K_0}(\tau), x)}= \zeta_{(\tau, x)}$, 
$\delta_{(H,  \sigma_{K_0}\tau))}^{ (\{1\}, \sigma_{K_0}( \phi))}
=\delta_{(H,  \tau)}^{ (\{1\},  \phi)}$ by Lemma \ref{isomoves} (iv) and  
$\langle 
\Res|_H\chi, \tau \rangle \zeta_{(\tau, x)} \in {\mathbb Z}$.

For $\chi \in \Irr (G, K)$,  $\phi \in \IBr(H,K)$, 
$$d_{G}^{(x,  e_Q)}(\chi)  = \sum_{\phi \in \IBr(H, e_Q, K)} 
\delta_{(G,  \chi)}^{ (x,  \phi)}\phi_0,$$

The compatibility condition (b) of Definition \ref{defnisot} 
is easily seen to follow from  the linearity of the maps $I^Q$, the claim 
and Lemma  \ref{isomoves}

\begin{eqnarray*} I^Q_{p'} \circ d_{G}^{(x,  e_Q)}(\chi) &=&  
\sum_{\phi \in \IBr(H, e_Q, K)} \delta_{(G,  \chi)}^{ (x,  \phi)}\sigma_{K_0} 
I^Q(\phi_0) \\
&=& \sum_{\phi \in \IBr(H, e_Q, K)} \delta_{(G,  \sigma_{K_0}(\chi))}^
{ (x,  \sigma_{K_0}
(\phi))} I^Q(\phi _0)\\
&=& \sum_{\phi \in \IBr(H, e_Q, K)} \delta_{(G,  \sigma_{K_0}(\chi))}^
{ (x,  \sigma_{K_0}
(\phi))} \sigma_{K_0}(\phi)\\
&=& \sum_{\phi' \in \IBr(H, \sigma_{K_0} (e_Q), K)} 
\delta_{(G,  \sigma_{K_0}(\chi))}^{ (x, \phi')}\phi'_0\\
&=& d_{G}^{(x,  \sigma_{K_0}(e_Q))}  I^{\{1\}}(\chi).
\end{eqnarray*}
 \end{prf1}

\

\begin{prf2}   Let $(K, \Oo, k)$ be a $p$-modular system such that  
$k$  is an algebraic closure of ${\mathbb F}_p$ and $ K$ contains a primitive
$|G|$-th root of unity.  By the theorem,  there is a perfect isometry between 
$\Oo G b $ and  $\Oo G \sigma_{K_0}(b) $. Hence, by 
\cite[Th\'eor\`eme 1.4]{broue-isome}, there is an $\Oo$-algebra isomorphism 
$ f : Z(\Oo  G \sigma_{K_0}(b)  )  \to Z(\Oo  G b)$. This 
induces a  $k$-algebra  isomorphism  $ f : Z(kG \sigma(\bar b))  
\to Z(kG \bar b) $.   
On the other hand, $\sigma $ induces a ring isomorphism 
$\sigma : Z(kG \bar b) \to Z(kG\sigma ( \bar b))$ such that  for all 
$ a\in  Z(kG \bar b)$ and all $\lambda \in k $, 
$\sigma (\lambda a) = \lambda^p \sigma(a)$. 
Thus, $\sigma \circ f :  Z(kG \sigma (b)) \to   Z(kG \sigma (b)) $ 
is a ring automorphism    which  satisfies  $\sigma \circ f (\lambda a ) = 
\lambda^p \sigma  \circ f   (a)  $ for all $ a\in  Z(kG \sigma ( \bar b)) $ 
and all  $\lambda \in k $.   By \cite[Lemma 2.1]{Ke1}, the fixed points 
$ (Z(kG\sigma ( \bar b)))^{\sigma  \circ f}  $
 of $ Z(kG\sigma ( \bar  b))  $  under $ \sigma  \circ f $ are an  
${\mathbb F}_p$-subspace of $Z(kG\sigma ( \bar b))$ such that   
$Z(kG\sigma( \bar b)) \cong  k \otimes_{{\mathbb F}_p}  
(Z(kG\sigma ( \bar b)))^{\sigma  \circ f}$ as $k$-vector spaces.
Since $ \sigma  \circ f$  is a homomorphism of rings, 
$ (Z(kG\sigma ( \bar b)))^{\sigma  \circ f}  $  is an ${\mathbb F}_p$-algebra. 
Thus
$Z(kG \sigma ( \bar b))$  and hence $Z(kGb)$   has an ${\mathbb F}_p$-form. 
\end{prf2}

\

\begin{prf3} Let $b$ be a block of $kG$ with defect $d$.   
Since ${\rm dim}_k(Z(kGb)) ={\rm dim}_K(Z(KGb))=|\Irr(G,b, K)|  $,  
by  \cite[Theorem 1]{brauerfeit} the  $k$-dimension of $Z(kGb)$ is bounded   
by  
$m$.   By   Corollary \ref{CPWcrl}, $Z(kGb)$  has  a $k$-basis  
such that the   
multiplicative  constants of $Z(kGb)$  with respect to this basis  are all  
in ${\mathbb F}_p$.  Thus there  are at most  $p^{m^3}$  possibilities   
for  the isomorphism  type of $Z(kGb)$.
\end{prf3}

\begin{rms}(i) 
\rm{The proof of Theorem \ref{mainthm} can be readily adapted   
to prove that  between  any  pair of Galois conjugate blocks  there is a  
global isotypy in the sense of  \cite[1.9]{Boltje/Xu:2008}. It is not known 
whether 
there is a $p$-permutation equivalence 
(cf. \cite[Definiton 1.3]{Boltje/Xu:2008},  
\cite[Definiton 1.3]{Linckelmann:p-perm})  between    any  pair of Galois  
conjugate blocks.} 

(ii) \rm {By Lemma  \ref{isomoves},  
the isometries  between various blocks in Theorem \ref{mainthm}   
all appear without  signs, which seems to render even  more surprising  the 
fact that  Galois conjugate blocks 
need not be Morita equivalent \cite{Benson/Kessar:2007}.}

\end{rms}

\end{document}